
\baselineskip=14pt
\parskip=10pt

\magnification=\magstephalf

\def\1{{\overline{1}}}
\def\2{{\overline{2}}}
\parindent=0pt
\overfullrule=0in

\def\frac#1#2{{#1 \over #2}}
\centerline
{\bf 
Experimenting with Discrete Dynamical Systems
}
\bigskip
\centerline
{\it George SPAHN and Doron ZEILBERGER}
\centerline
\qquad \qquad \qquad 

{\it Dedicated to Saber Elaydi  on his  $80^{th}$ birthday, and to Gerry Ladas on his $85^{th}$  birthday}

\bigskip

{\bf Abstract}:  We demonstrate the power of Experimental Mathematics and Symbolic Computation to study intriguing problems on rational difference equations, studied 
extensively by Difference Equations giants, Saber Elaydi and Gerry Ladas (and their students and collaborators). In particular we rigorously prove some fascinating conjectures made by Amal Amleh and
Gerry Ladas back in 2000. For other conjectures we are content with semi-rigorous proofs. We also extend  the work of Emilie Purvine (formerly Hogan) and
Doron Zeilberger for rigorously and semi-rigorously proving global asymptotic stability  of arbitrary rational difference equations (with positive coefficients), and
more generally rational transformations of the positive orthant of $R^k$ into itself.

{\bf Keywords}:  Non-Linear Difference Equations, Discrete Dynamical Systems, Semi-Rigorous Proofs

{\bf 2020 MSC: 39A23, 39A30}

{\bf Preface}: Difference Equations, linear, but especially non-linear, are not only so useful in mathematical biology and elsewhere, but
are fascinating to study for their own sake. Recall that already the {\bf first order} {\it logistic} difference equation
$$
x_{n+1} = \lambda x_{n}(1-x_n) \quad, \quad 0 \leq \lambda \leq 4 \quad,
$$
introduced by Sir Robert May, lead to {\it chaos}, one of the central paradigms of our time, as well as  to {\it period-doubling}, and the {\it Feigenbaum constants}.
Here the action takes place  in the finite interval $0 \leq x \leq 1$.

Often in population dynamics one encounters higher-order difference equation where $x_{n+1}$ is a {\bf rational} function of the previous $k$ values.
$$
x_{n+1} \, = \, \frac{a_0+ \sum_{i=1}^{k} a_i x_{n+1-i}}{b_0+ \sum_{i=1}^{k} b_i x_{n+1-i}} \quad,
$$
where the coefficients $a_0,a_1, \dots, a_k$ and $b_0,b_1, \dots, b_k$ are assumed {\bf non-negative} (and of course at least one of them is strictly positive at the bottom),
and with positive {\bf initial conditions} $x_1, \dots , x_k$. Here the action takes place in the infinite interval $0 < x <\infty$.

These difference equations have been studied extensively, with a deep theory, by the two `birthday boys' and their many disciples, see [AL],[CL], [E1], [E2], [KoL], [KuL], and
references thereof.

Twelve years ago one of us (DZ), in collaboration with his then PhD student, Emilie Hogan (now Purvine), initiated the use of {\it symbolic computation},
and computer-generated (rigorous!) proofs to study such difference equations, but only those whose solutions always converge to a {\bf unique}
{\it equilibrium point}, i.e. for which there exists a unique $\bar{x}>0$ such that for {\it any} (positive) {\bf initial conditions},
$x_1, \dots, x_k$, one has $\lim_{n \rightarrow \infty} x_n \, = \, \bar{x}$, and the challenge was to prove it rigorously.
Note that to prove that a candidate fixed point $\bar{x}$ is {\it locally stable} can be easily done using standard techniques (see below).

In [AL] several fascinating conjectures were made (with monetary prizes offered by Gerry Ladas for some of them, alas with the very ungenerous
deadline of Jan. 1, 2002, long passed). Here is the first one:

{\bf Conjecture ([AL])}: Prove that every positive solution of the difference equation
$$
x_{n+1} \, = \, \frac{x_{n-1}}{x_{n-1}+x_{n-2}} \quad,
$$
converges to a period two solution of that equation of the form
$$
\dots, \phi, 1-\phi, \dots \quad,
$$
with $0 \leq \phi \leq 1$.

In this paper we extend, and improve, the method of [HZ], and then extend this methodology to  rigorously 
prove (with the aid of our beloved silicon servants) some of the fascinating conjectures in [AL].
For other ones we will be content with {\it semi-rigorous} proofs (see [Z] for the concept). We will argue
that often such proofs suffice, since we know that there {\bf exists} a rigorous proof (or disproof), but
it would be a waste of the computer's time and memory to find it, since the probability that the
semi-rigorous proof was a false positive is negligible.

{\bf Accompanying Maple packages}

This article is accompanied by two Maple packages.

$\bullet$ {\tt DRDS.txt}, to experiment, numerically, and symbolically, with solutions of rational difference equations of {\it any} order, and
for proving, if possible, {\it rigorously}, but more often (if the order is three and up) {\it semi-rigorously}, global asymptotic stability. This latter part is
a continuation, and improvement, of the pioneering work in [HZ], that focused on second-order difference equations.

$\bullet$ {\tt AmalGerry.txt}, to prove rigorously (and in some cases, {\it semi-rigorously}), some of the intriguing and tantalizing conjectures
made by Amal Amleh and Gerry Ladas in 2001 [AL].

Both packages, and numerous input and outputs files, are viewable (and downloadable!) from the front of this article

{\tt https://sites.math.rutgers.edu/\~{}zeilberg/mamarim/mamarimhtml/dds.html } \quad .

{\bf Experimenting with Some Random Rational Difference Equations}

One of the purposes of this article is to serve as a {\it tutorial} on our Maple package {\tt DRDS.txt}, that in addition to {\it proving}
can be also used as a {\it calculator} for exploring and experimenting, {\it numerically}, with rational difference equations.

To get a feel of how a {\it typical} rational difference equation looks like, use  procedure {\tt RRDE(x,n,k,d,A)},
where the {\bf input parameters} are:

$\bullet$ {\tt x} and {\tt n} are {\it symbols} that correspond to the notation $x_n$ in `humanese'.

$\bullet$ {\tt k} is the {\it order} of the difference equation

$\bullet$ {\tt d} is the {\it degree} of the numerator and denominator (in this paper we will focus on the degree one case, but the Maple package can handle any degree).

$\bullet$ {\tt A} is a positive integer.

The {\bf output} is such a random rational difference equation whose numerator and denominators have integer coefficients drawn from $\{1, \dots, A\}$.

For example, typing

{\tt F:=RRDE(x,n,3,2,30);} \quad, 

may give you something like this (of course, every time you get a different difference equation, since this is random)
$$
x_{n+1}=
\frac{17 x_{n}^{2}+2 x_{n} x_{n -2}+20 x_{n} x_{n -1}+21 x_{n -2}^{2}+25 x_{n -2} x_{n -1}+17 x_{n -1}^{2}+8 x_{n}+23 x_{n -2}+6 x_{n -1}+13}{4 x_{n}^{2}+4 x_{n} x_{n -2}+29 x_{n} x_{n -1}+10 x_{n -2}^{2}+4 x_{n -2} x_{n -1}+4 x_{n -1}^{2}+9 x_{n}+6 x_{n -2}+19 x_{n -1}+9}
\, .
\eqno(1)
$$

To get the first $N+k$ terms of such a difference equation with given initial conditions, type:

{\tt L:=OrbD(F,x,n,INI,N);}

where  {\tt F} is the difference equation in the above format $x_{n+1}=F$, {\tt x} and {\tt n} are the symbols used to express it, {\tt INI} is a list of 
of length $k$, and {\tt N} is
the number of extra terms (so the output list is of length $N+k$). For example, typing:

{\tt OrbD(x[n]+x[n-1],x,n,[1,1],10);} 

would yield

{\tt [1, 1, 2, 3, 5, 8, 13, 21, 34, 55, 89, 144]}   \quad.

Going back to the above complicated (random) third-order difference equation denoted by {\tt F} above, typing

{\tt L:=OrbD(F,x,n,evalf([21,27,39]),1000):}

would give you in {\it floating-point}, the first $1003$ terms of the sequence with initial conditions $x_1=21,x_2=27,x_3=39$. Note the colon {\tt :}, as opposed to the semi-colon, {\tt ;}, since
we really don't want to see all of them, we are only interested in the {\it long-term behavior}, i.e. whether {\bf ultimately} it converges to one number (or in the long-run to a periodic solution (see below)).

Typing {\tt L[-1];} and {\tt L[-2];} would give respectively, the $1002^{th}$ and $1003^{th}$ terms:
$$
1.6358881124186260402\dots \quad , \quad 1.6358881124186260402\dots \quad,
$$
indicating that they are very close, hence, at least with the above initial conditions, the sequence seems to have  a limit. Experimenting with randomly chosen (positive) initial
conditions, we see again and again, that it seems to converge to the same number. Hence, completely based on {\it numerics}, we are safe to make the following conjecture.

{\bf Random Conjecture}: For {\it any } positive initial conditions $x_1,x_2,x_3$, the terms of the sequence satisfying the difference equation $(1)$ converge to a certain
algebraic number (that can be easily found), whose floating-point value  is $1.6358881124186260402\quad$.

So far we only used {\it numerical computations}. Later  on we will use {\it symbolic computation} to actually prove it, either rigorously or semi-rigorously.

Experimenting with many other random cases, and several initial conditions, you get again and again this phenomenon that the sequences seem to converge to the
{\bf same} number, let's call it $\bar{x}$, i.e. that it is a {\bf global equilibrium}. To find out its value, one replaces $x_{n+1},x_n,x_{n-1}, \dots, x_{n-k+1}$ by $\bar{x}$
getting a one-variable {\bf polynomial} equation in $\bar{x}$, asking Maple to {\tt solve} it, and if all goes well only getting {\bf one} real and positive solution
(of course, the resulting equation would have $d+1$ complex roots).

Let's pick  another random example, this time a {\it third-order} difference equation, and to make it simpler, let's have the numerator and denominator of degree $1$. Typing

{\tt T:=RRDE(x,n,3,1,30);} 

yielded (this time)
$$
x_{n+1} \, = \,\frac{17+5 x_{n}+24 x_{n -1}+16 x_{n -2}}{23+4 x_{n}+19 x_{n -1}+2 x_{n -2}} \quad .
\eqno(2)
$$

Taking random initial conditions (in this case $x_1=11,x_2=27,x_3=37$, and typing

{\tt L:=OrbD(T,x,n,evalf([11,27,37]),1000):  L[-1],L[-2];} 

gives:

$$
1.37466571564383380885297479\dots \quad , \quad 1.37466571564383380885297479\dots \quad,
$$
and similarly for  many other randomly chosen initial conditions. To find the {\bf exact} value of this (so far conjectured) equilibrium point,
solve
$$
\bar{x} \, = \,\frac{17+5 \bar{x}+24 \bar{x}+16 \bar{x}}{23+4 \bar{x} +19 \bar{x} +2 \bar{x}} \quad,
$$
getting
$$
\bar{x}= \frac{17+45 \bar{x}}{23+25\bar{x}} \quad,
$$
that simplifies to the {\bf quadratic} equation:
$$
25 \bar{x}^{2}-22 \bar{x} -17 \, = \, 0 \quad , 
$$
whose roots are:
$$
\left[\frac{11}{25}+\frac{\sqrt{546}}{25}, \frac{11}{25}-\frac{\sqrt{546}}{25}\right] \quad,
$$
that in decimals are:
$$
[1.374665716, -0.4946657156] \quad,
$$

discarding the negative root, we got the exact value, namely  $\frac{11}{25}+\frac{\sqrt{546}}{25}$.

{\bf The Amleh-Ladas Fascinating conjectures}

Procedure {\tt RRDE} artificially made all coefficients strictly positive, and hence it turns out that, {\it generically}, one gets rather boring limiting
behavior, i.e. convergence to a unique  positive equilibrium point, or phrased otherwise, a {\it limiting period-one solution}.

In [AL], Amal Amleh and Gerry Ladas made the following intriguing conjectures, that exhibited far more interesting long-term behavior.

{\bf Conjecture 1}: For any positive {\it initial conditions} $x_1,x_2,x_3$, 
$$
x_{n+1}= \frac{x_{n-1}}{x_{n-1}+x_{n-2}} \quad, \quad n \geq 3 \quad,
$$
the sequence $\{x_n\}$ converges to a {\bf period-two} solution of the form
$$
\dots \,, \, \phi \, , \, 1- \phi \, , \, \dots \quad,
$$
with $0 \leq \phi \leq 1$.

{\bf Conjecture 2}: For any positive {\it initial conditions} $x_1,x_2,x_3$, 
$$
x_{n+1}= \frac{x_{n}+x_{n-2}}{x_{n-1}} \quad, \quad n \geq 3 \quad,
$$
the sequence $\{x_n\}$ converges to a {\bf period-four} solution of the form
$$
\dots \,, \, \phi \, , \, \psi \, , \, \frac{\phi+\psi^2}{\phi \psi -1} \, , \, \frac{\phi^2+ \psi}{\phi\psi -1} \, \dots \quad ,
$$
with $\phi,\psi \in (0, \infty)$, and $\phi\psi>1$.

{\bf Conjecture 3}: For any positive {\it initial conditions} $x_1,x_2,x_3$, 
$$
x_{n+1}= \frac{1+x_{n-2}}{x_{n}} \quad, \quad n \geq 3 \quad,
$$
the sequence $\{x_n\}$ converges to a {\bf period-five} solution of the form
$$
\dots \,, \, \phi \, , \, \psi \, , \, \frac{1+\phi}{\phi \psi -1} \, , \, \phi\psi-1 \, , \,  \frac{1+ \psi}{\phi\psi -1} \, \dots
$$
with $\phi,\psi \in (0, \infty)$, and $\phi\psi>1$.

{\bf Conjecture 4}: For any positive {\it initial conditions} $x_1,x_2,x_3$, 
$$
x_{n+1}= \frac{1+x_{n}}{x_{n-1}+x_{n-2}} \quad, \quad n \geq 3 \quad
$$
The sequence $\{x_n\}$ converges to a {\bf period-six} solution of the form
$$
\dots \,, \, \phi \, , \, \psi \, , \, \frac{\phi}{\psi} \, , \, \frac{1}{\phi} \,, \, \frac{1}{\psi} \, , \, \frac{\phi}{\psi} \, , \dots \quad,
$$
with $\phi,\psi \in (0, \infty)$.

We will later show how to prove these using {\bf symbolic computation}, but for now, let's confirm them numerically, using procedure {\tt OrbD}.

These difference equations  are hard-coded in procedure {\tt LadadDB(x,n)}, that contains $15$ interesting difference equations, the first four, namely

{\tt LadasDB(x,n)[1] \quad, \quad LadasDB(x,n)[2] \quad, \quad LadasDB(x,n)[3] \quad, \quad LadasDB(x,n)[4]}, 

correspond to the difference equations featured in the above four conjectures.

For example, typing

{\tt T:=LadasDB(x,n)[1]: evalf(OrbD(T,x,n,evalf([1,2,3]),1000)[-3..-1],10);}

gives
$$
[0.7012220196, 0.2987779809, 0.7012220196] \quad,
$$

while

{\tt T:=LadasDB(x,n)[1]: evalf(OrbD(T,x,n,evalf([11,25,34]),1000)[-3..-1],10);}

yields
$$
[0.9348089961, 0.06519100396, 0.9348089961] \quad .
$$

Readers are welcome to experiment with many other initial conditions, and similarly for the other difference equations featured in Conjectures $2$,$3$, and $4$, to numerically (empirically) confirm these
intriguing conjectures. Of course this `only' gives numerical confirmation (that the authors of [AL] must have already done, but probably using numerical software rather than Maple).
We will later see how to prove them either rigorously or semi-rigorously.

{\bf Recalling some basics and More Numerical Explorations}

While it is highly non-trivial, in general, to prove that {\it every} choice of initial conditions will make the solution sequence converge to an equilibrium, it
is purely routine, {\it today}, to decide whether it is true when the initial conditions are not too far from that equilibrium.

The first step (already recalled in [HZ]) is to convert a $k$-th order difference equation in $(0, \infty)$ to a {\bf first}-order difference equation in $(0,\infty)^k$.

The $k$-th order difference equation
$$
x_{n+1}= F(x_n,x_{n-1}, \dots, x_{n-k+1}) 
$$
becomes the transformation
$$
\left [ \matrix{x_1 \cr x_2 \cr \dots \cr x_k} \right ] \rightarrow
\left [ \matrix{x_2 \cr x_3 \cr \dots \cr x_k \cr F(x_k,x_{k-1},x_{k-2}, \dots, x_1) } \right ]  \quad.
$$
From this point of view, the more general problem is to investigate whether, given a general {\bf rational} transformation from $(0,\infty)^k$ into $(0,\infty)^k$  of the form
$$
\left [ \matrix{x_1 \cr x_2 \cr \dots \cr x_k} \right ]\rightarrow 
\left [ \matrix{R_1(x_1, \dots, x_k)  \cr R_2(x_1, \dots, x_k) \cr \dots \cr R_k(x_1,\dots, x_k)} \right ] \quad,
$$
where $R_1(x_1, \dots, x_k) , \dots , R_k(x_1, \dots, x_k)$ are {\bf rational functions} of their arguments.
In order to prove that $\bar{x}$ is the limit of every solution of the original difference equation, one has to prove that the {\bf orbit}, starting at any point in $(0,\infty)^k$ of the
resulting transformation (as constructed above) converges to $(\bar{x}, \bar{x}, \dots, \bar{x})$.

Procedure {\tt Targem} in the Maple package {\tt DRDS.txt} converts a difference equation to a transformation.

Recall that for any transformation ${\bf x} \rightarrow {\bf F}({\bf x})$ in $R^k$, a point ${\bf x}_0$ is called an {\it equilibrium point} if it is a {\bf fixed point}.
$$
{\bf F}({\bf x}_0) = {\bf x}_0 \quad .
$$

If the transformation is, like in our case, rational, this gives a system of $k$ polynomial equations with $k$ unknowns, that at least, {\it in principle}, but often also in practice (for small $k$)
is fully solvable (in the realm of algebraic numbers, if all the coefficients are integers).

In order to generate random examples to experiment with, use procedure

{\tt RRT(x,k,d,A)},

where {\tt x} is a symbol, {\tt k} is the {\bf dimension} of the space, {\tt d} is the degree of the numerator and the denominator, and A is a positive integer, such that
the coefficients are randomly chosen from $\{1,2, \dots, A\}$. For example,

{\tt T:=RRT(x,2,1,30);} \quad,

might give 
$$
\left[\frac{18+20 x_{1}+24 x_{2}}{11+19 x_{1}+25 x_{2}}, 
\frac{26+29 x_{1}+28 x_{2}}{29+14 x_{1}+18 x_{2}}\right] \quad .
$$

To get {\bf all} the fixed points (including complex and with negative coordinates), type

{\tt FP(T,x);}

In this example, you would get a big mess, so let's take the floating-point version, and type

{\tt evalf( FP(T,x),10);}

getting

{\tt [[0.5983411214, -1.834086341], [1.100408318, 1.394961226], [-0.9483476940, 0.159782893], [-100.8951386, 76.30710943]]} \quad.

We are only interested in points in $(0,\infty)^2$, so the only point that we are interested in is the second one. To get such points
with all positive coordinates right away, type

{\tt evalf( FPp(T,x),10);}

getting, indeed,

$$
\{[1.100408318\dots, 1.394961226\dots] \} \quad .
$$

We next ask whether it is {\it locally stable}. As is well-known, and fairly easy to see (e.g. [KL]), one computes the {\bf Jacobian matrix},
then for each of the candidate points (those with positive coordinates)  {\bf plugs-it-in}
then computes the eigenvalues of this numerical matrix, and if all of them have absolute value less than one, then we 
know {\bf for sure} that the examined equilibrium point  is {\it locally stable}. This is an important first step before we can prove
{\bf global stability}, since every globally stable fixed point must be, first of all, a local one.

Of course if there are more than one locally stable fixed points in $(0,\infty)^k$ there is no hope that any of them is global, but
it is still nice to know all of them.

This is all done, thanks to Maple, {\it automatically}, with procedure {\tt LSFP(T,x)}.

For example, with the above {\tt T}, we would get that this point is indeed locally stable. It also gives you
the (floating-point) versions of the eigenvalues (just for the record).

So typing

{\tt evalf(LSFP(T,x),10);}

gives
$$
\{[[ 1.100408318,  1.394961226], [
 0.01399544579+ 0.07999899783 \,i, 
 0.01399544579- 0.07999899783 \,i]]\} \quad,
$$

indicating that indeed there is only one locally stable fixed point of our transformation, and also giving the eigenvalues.

Before trying to prove {\it rigorously} (that takes lots of effort!), or even {\it semi-rigorously} (that also takes some effort, see below) it is very
easy to conjecture whether there (most probably) is a globally stable fixed point, and to actually find it. In fact, we don't need to find the locally stable
fixed points. The completely numeric and {\it empirical} procedure

{\tt CoGSFP(T,x,K1,K2)} \quad,

takes as input a transformation {\tt T} and picks {\tt K2} random points in $(0,\infty)^k$, and for each of them computes the orbit of length {\tt K1}.
If for each of these orbits the difference between the last two terms is very small, and they all give the same point, we know
{\bf non-rigorously}, but  (almost) {\bf certainly} that there is a globally stable fixed point, and we know its floating point version.
To get its exact value, as an algebraic number, you need {\bf FPp} mentioned above.

{\bf The Proof Strategy}

This is essentially what was done in [HZ], but with a new implementation, and with a {\it semi-rigorous} option.

Suppose that we have a (rational) transformation, $T$,  from $(0,\infty)^k$ to  $(0,\infty)^k$ , and that we already have ample empirical/numerical evidence that
a certain {\it candidate point}, (gotten either from {\tt LSFP} or {\tt CoGSFP}), $\bar{x}$, is a globally stable fixed point of $T$ in $(0,\infty)^k$.

Our goal in life is to prove that for {\it every} {\bf initial point}, ${\bf x}_0 \in (0,\infty)^k$,

$$
\lim_{n \rightarrow \infty} T^n ({\bf x}_0) \, = \, \bar{x} \quad.
$$

Let $|x|$ be the usual Euclidean norm. The above statement is equivalent to

$$
\lim_{n \rightarrow \infty} |T^n ({\bf x}_0)-  \bar{x}|^2 \, = \, 0 \quad.
$$

Suppose that we can come-up with a {\bf positive} integer $r$, and some real $\alpha>1$, such that, for any point ${\bf x} \in (0,\infty)^r$,  we have the
inequality
$$
\alpha \,|T^r ({\bf x})-  \bar{x}|^2 \leq |{\bf x}-  \bar{x}|^2 \quad,
$$
then we know that this sequence of `distance-squared from the fixed point' shrinks (by at least the factor $1/\alpha$), every $r$-th iteration compared to what it was. This would automatically entail what we want.
Note that this is only a {\it sufficient} condition, and there is no {\it a priori} reason (that we know of, at least), that such a real $\alpha>1$ and an integer $r \geq 1$ exist, but if we are
{\bf lucky} enough to find a candidate, and then, succeed in proving it, we are be done!

Again, we {\bf first} investigate things {\it numerically}.
We start with $r=1$, and then see whether for many initial points, the resulting orbit has the property that
the distance-squared from $\bar{x}$ shrinks every $r$-th iteration. Once we get a successful candidate $r \geq 1$, we have to {\bf prove} it.
For the sake of definiteness, and not to clutter the Maple code with another parameter, we decided to take $\alpha=\frac{101}{100}$.

So we have to prove (either rigorously or semi-rigorously) that for any ${\bf x} \in (0,\infty)^k$, we have the {\bf inequality}

$$
\alpha |T^r({\bf x})-{\bar x}|^2  \leq |x -\bar{x}|^2 \quad.
$$

in other words
$$
|x -\bar{x}|^2 - \alpha|T^r({\bf x})-{\bar x}|^2 \geq 0 \quad.
$$

Now the left-side is a (usually rather complicated) {\bf rational function} of $x_1, \dots, x_k$, (recall that $x= (x_1, \dots, x_k)$).
Simplifying, we get a {\bf denominator} that is
a {\bf perfect square}, and hence automatically positive. The numerator is a certain (usually complicated) polynomial, and everything boils down to proving that
this polynomial, let's call it $P(x_1, \dots, x_k)$, is non-negative. In other words we need to {\bf minimize} $P$ in the region $(0,\infty)^k$
and prove that it is $\geq 0$ (in fact if it is $\geq 0$ it must be $0$, since the value of the above left side is $0$ when $x=\bar{x}$).  

But this is a {\bf routine} multivariable calculus exercise!, that Maple (and Mathematica, and Sage),  know how to do.
Alas, in applications to our problems, this polynomial turns out to be too complicated for $k >2$, so we can rigorously prove 
global stability for $k=2$, but it would take too long (on our modest laptops) to do it for $k \geq 3$. So we opt to do it numerically, checking it for many random points, and
since we know that there {\bf exists} a way, at least in principle, to prove it rigorously, why bother? This {\it semi-rigorous}
approach to mathematical deduction was proposed by one of us thirty years ago [Z].

To get a fully rigorous proof, use procedure {\tt GSFPv(T,x,K);}, where {\tt K} is the maximum $r$ we are willing to take.
That works well with two dimensions, but for higher dimensions, it takes way too long.
One should use instead {\tt GSFPvSR(T,x,K);}, in order to get a semi-rigorous proof. If none is found, it returns {\tt FAIL}.

{\bf Sample output for proving (Rigorously and Semi-Rigorously Global Stability)}

$\bullet$ If you want to see $20$  theorems that state that certain second-order difference equations  always converge to the unique stable
equilibrium, with {\bf fully rigorous} proofs, look here:

{\tt https://sites.math.rutgers.edu/\~{}zeilberg/tokhniot/oDRDS1.txt} \quad .

$\bullet$ If you want to see $10$  theorems that state that certain third-order difference equations  always converge to the unique stable
equilibrium, with {\bf semi-rigorous} proofs, look here:

{\tt https://sites.math.rutgers.edu/\~{}zeilberg/tokhniot/oDRDS2.txt} \quad .

$\bullet$ If you want to see $5$  theorems that state that certain fourth-order difference equations  always converge to the unique stable
equilibrium, with {\bf semi-rigorous} proofs, look here:

{\tt https://sites.math.rutgers.edu/\~{}zeilberg/tokhniot/oDRDS3.txt} \quad .

(Warning: the file is large!)

{\bf The Amleh-Ladas Conjectures}

We now apply the proof strategy described above to the fascinating Amleh-Ladas conjectures. Recall conjecture 1:

{\bf Conjecture 1}: For any positive {\it initial conditions} $x_1,x_2,x_3$, 
$$
x_{n+1}= \frac{x_{n-1}}{x_{n-1}+x_{n-2}} \quad, \quad n \geq 3 \quad,
$$
the sequence $\{x_n\}$ converges to a {\bf period-two} solution of the form
$$
\dots \,, \, \phi \, , \, 1- \phi \, , \, \dots \quad,
$$
with $0 \leq \phi \leq 1$.

To tackle this problem, we first view the sequence as a map from $R^3$ to $R^3$. $$T:(x,y,z) \rightarrow (y,z,\frac{y}{x+y})$$ We hope to show that for any start point $(x_0,y_0,z_0)$ with $x_0,y_0,z_0 > 0$, this dynamical system, $T$ converges to a point somewhere on the line segment parametrized by $(t,1-t,t), t\in [0,1]$.

To measure how close points are to this line, we use the norm $$ v(x,y,z) = 1+x^2+y^2+z^2-x-2y-z+xy-xz+yz $$ \quad .
This is the square of the Euclidean distance from $(x,y,z)$ to the line. To show that $T$ converges to a point with $v$-norm equal to 0, we consider the objective function $$ F(x,y,z) = v(x,y,z) - v(T(x,y,z)) $$
If we can show that this objective function is always $\geq 0$ for all positive $(x,y,z)$, then we conclude that applying $T$ cannot increase the $v$-norm. Unfortunately this objective function is sometimes negative for this particular $T$ and $v$. To fix it, we replace $T$ with $T^3$, three consecutive applications of $T$. We also apply $T$ to both points for smoothing. This gives $$ F(x,y,z) = v(T(x,y,z)) - v(T^4(x,y,z)) $$ \quad .

It is now time to set maple to work and show that the objective function is nonnegative. We use the function NLPSolve from the Optimization package, to minimize the objective function. This function requires an initial point, and uses iterative methods to search for improvements on the initial point using floating point arithmetic. The maple documentation for this function is given here:

{\tt https://www.maplesoft.com/support/help/maple/view.aspx?path=Optimization\%2FNLPSolve }

We ran the solver with 36 different initial points and each time it returned that the minimum was within an acceptable range $(\pm 10^{-6})$ of 0. This completes a semi-rigorous proof of Conjecture 1. To run our code yourself, download the maple package 
{\tt AmalGerry.txt}, and type:

{\tt run\_nlp(T4,v4,1,4,3);} \quad .

The first argument asks for the transformation on $R^3$. T4 is the transformation T defined above; the 4 is a reference to the fact that this transformation was equation 4 on the Amleh-Ladas paper. The second argument is the norm to be used; v4 corresponds to T4. The third and fourth arguments specify how many iterations of T should be applied to the initial point when creating the objective function. The last argument determines the amount of different initial points that are tested.

We used this method to semi-rigorously prove conjectures 1 through 4. Take a look at the maple package {\tt AmalGerry.txt} for more details!

\vfill\eject

{\bf Using Maple's symbolic minimizer}

The NLPSolve function is not a mathematical proof of minimization. Perhaps it could be made into one by taking a very close look at the algorithm in the code, but this would be very tedious. In this section we attempt to use maple's built-in symbolic minimizer to produce a fully rigorous proof of the minimization. The main limitation here is computational resources, so we experiment with different norms and parameters.

Still looking at the $T$ and $v$ from the previous section, maple is not able to symbolically compute the minimum of  $$ F(x,y,z) = v(T(x,y,z)) - v(T^4(x,y,z)) $$ in a reasonable amount of time.
Instead we define simpler norms $$v_{xy}= (x+y-1)^2$$ $$ v_{yz}=(y+z-1)^2$$
Let $$F_{xy}=v_{xy}(x,y,z)-v_{xy}(T^4(x,y,z))$$
The denominator of $F_{xy}$ turns out to be $$(xz + yz + y)^2(y + z)^2$$ which is never negative! Thus we just instruct maple to minimize the numerator, which is a degree 8 polynomial in the variables $x,y,z$. Maple symbolically computes that the minimum of this polynomial is 0, so we have a rigorous proof that $T$ converges to a point with $v_{xy}$ norm equal to 0. The same process works for $v_{yz}$. The only points where $v_{xy}$ and $v_{yz}$ are both 0 is exactly the line $(t,1-t,t)$. This completes a rigorous mathematical proof of Conjecture 1. The authors attempted to use a similar approach for the other conjectures however we lacked the computational resources to execute the code.

{\bf Periodic Difference Equations}

As mentioned in [AL] (p. 71), and [KL] (p. 628), the following (second-order) {\bf Lynnes difference Equation}

$$
x_{n+1} \, = \, \frac{1+x_n}{x_{n-1}} \quad
$$

{\bf always} has period five, regardless of the initial conditions. This is easily confirmed with our Maple package:

Entering:

{\tt OrbD((1+x[n])/x[n-1],x,n,[a,b],5);}

immediately gives
$$
\left[a, b, \frac{1+b}{a}, \frac{a +1+b}{a b}, \frac{a +1}{b}, a, b \right] \quad,
$$

meaning that for {\bf symbolic}, i.e. {\it general}, initial conditions, things get repeated every five iterations.
Of course this is easy enough to do {\it by hand}.

Also mentioned in [KL] (p. 628, Eq. (25) there), that the following (third-order) difference equation

$$
x_{n+1} \, = \, \frac{1+x_n+x_{n-1}}{x_{n-2}} \quad ,
$$

is always of period eight, regardless of the initial conditions. Indeed, typing:

{\tt OrbD((1+x[n]+x[n-1])/x[n-2],x,n,[a,b,c],8);}

yields
$$
\left[a, b, c, \frac{1+c +b}{a}, \frac{c a +a +b +c +1}{a b}, 
\frac{a b +c a +b^{2}+b c +a +2 b +c +1}{a b c}, 
\frac{c a +a +b +c +1}{b c}, \frac{a +1+b}{c}, a, b, c\right] \, .
$$

This gave us the hope that the fourth-order difference equation
$$
x_{n+1} \, = \, \frac{1+x_n +x_{n-1}+x_{n-2}}{x_{n-3}} \quad ,
$$
is perhaps periodic? Alas, entering

{\tt
L:=OrbD((1+x[n]+x[n-1]+x[n-2])/x[n-3],x,n,[1,1,1,1],1000): member(1,{op(5..nops(L),L)});
}

gives {\tt false}, meaning, that if there is a period, it would be larger than $1000$, so this difference equation is unlikely to be periodic.

It would be very interesting to discover such rational difference equations with higher periods, that do not
trivially follow from the known ones by `merging'.

{\bf Automated Discovery of Invariants that Imply that Every Solution is Bounded}

In the Kulenovic-Ladas fascinating book [KL], it is mentioned that the {\bf generalized Lynnes Equation}

$$
x_{n+1} \, = \, \frac{p+x_n}{x_{n-1}} \quad ,
$$

for any positive $p$ has the following invariant:
$$
I_n= (p+x_{n-1}+x_n) \left (1+ \frac{1}{x_{n-1}} \right ) \left (1+ \frac{1}{x_{n}} \right) \quad.
$$
We could have found it, {\it ab initio}, using procedure {\tt FindInv}:

{\tt FindInv((p+x[n])/x[n-1],x,n,2,3);}

We were able to find invariants for the higher-order difference equations
$$
x_{n+1}=\frac{p+x_n+x_{n-1}+ \dots+x_{n-k+2}}{x_{n-k+1}} \quad,
$$
for $k \leq 8$. They do get more and more complicated, see the output file

{\tt https://sites.math.rutgers.edu/\~{}zeilberg/tokhniot/oDRDS6.txt} \quad .

They all turn out to have positive coefficients. As noted in [KL] for the generalized Lynnes equation, but is also true in general, the 
existence of an invariant of the form
$$
\frac{P(x_n,x_{n-1}, \dots, x_{n-k+1})}{x_n x_{n-1} \dots x_{n-k+1}} \, = \, Constant \quad,
$$
with the coefficients of $P$ all positive
(all the ones we found were of that form) immediately implies that for any positive initial conditions, the solution sequence is always {\bf bounded}.

In fact, for the generalized Lynnes equation $x_{n+1}=\frac{p+x_{n}}{x_{n-1}}$, one can use {\bf discriminants} to predict, {\it a priori}, these
lower and upper bounds, for any positive $p$ and any positive initial conditions $x_1=a_1, x_2=a_2$. See the output file

{\tt https://sites.math.rutgers.edu/\~{}zeilberg/tokhniot/oDRDS4.txt} \quad .

For the analogous third-order  difference equation $x_{n+1}=\frac{p+x_{n}+x_{n-1}}{x_{n-2}}$, see:

{\tt https://sites.math.rutgers.edu/\~{}zeilberg/tokhniot/oDRDS5.txt} \quad .

{\bf Conclusion}

Using the great human-generated theory developed by Saber Elaydi, Gerry Ladas, and their many collaborators and disciples, and bringing into
the game both {\it numeric} and {\it symbolic} computation, we hope that we demonstrated the power of computer-kind to extend
the human efforts. Alas, even computers have their limits, and we advocate that often a {\it semi-rigorous} proof suffices, as
first preached in 1993 in [Z].

{\bf References}

[AL] A.M. Amleh and G. Ladas, {\it Convergence to Periodic Solutions}, 
Journal of Difference Equations and Applications {\bf 7}(2001), 621-631. \hfill\break

[CL] E. Camouzis and G. Ladas, {\it ``Dynamics of Third Order Rational Difference Equations''}, Chapman and
Hall/CRC press (2008).

[E1] Saber Elaydi, {\it ``Introduction to Difference Equations''}, Springer (2000).

[E2] Saber Elaydi, {\it Global dynamics of discrete dynamical systems and difference equations},
Difference equations, discrete dynamical systems and applications, 51-81,
Springer Proc. Math. Stat. {\bf 287}, Springer, 2019.

[HZ] Emilie Hogan (now Purvine) and Doron Zeilberger, {\it A New Algorithm for Proving Global Asymptotic Stability of Rational Difference Equations},
Journal of Difference Equations and Applications {\bf 18}(2012), 1853-1873. \hfill\break
{\tt https://sites.math.rutgers.edu/\~{}zeilberg/mamarim/mamarimhtml/gas.html} \quad.

[KoL] V. Kocic and G. Ladas, {\it ``Global Behavior of Nonlinear Difference Equations of Higher Order with Applications''}, Kluwer Academic Publishers (1993).

[KuL] M. R. S. Kulenovic and G. Ladas, {\it ``Dynamics of Second Order Rational Difference Equations}, Chapman and Hall/CRC press (2001).

[Z] Doron Zeilberger, {\it Theorems for a price: tomorrow's semi-rigorous mathematical culture},
Notices of the Amer. Math. Soc. {\bf 40}, 978-981. Reprinted in:
Math. Intell. {\bf 16} (1994) 11-14. \hfill\break
{\tt https://sites.math.rutgers.edu/\~{}zeilberg/mamarim/mamarimhtml/priced.html} \quad .

\bigskip
\hrule
\bigskip
George Spahn and Doron Zeilberger, Department of Mathematics, Rutgers University (New Brunswick), Hill Center-Busch Campus, 110 Frelinghuysen
Rd., Piscataway, NJ 08854-8019, USA. \hfill\break
Email: {\tt  gs828 at math dot rutgers dot edu} \quad, \quad {\tt DoronZeil] at gmail dot com}   \quad .

Written: {\bf June 21, 2023}. 

\end